\magnification=\magstep1
\hsize=16.5 true cm 
\vsize=25 true cm
\font\bff=cmbx10 scaled \magstep1
\font\bfff=cmbx10 scaled \magstep2
\font\smc=cmcsc10
\parindent0cm
\overfullrule=0pt
\font\boldmas=msbm10           
\def\Bbb#1{\hbox{\boldmas #1}} 

\centerline{\bfff Increasing singular functions with arbitrary} 
\smallskip
\centerline{\bfff positive derivatives at densely lying points}
\bigskip
\centerline{\bff Gerald Kuba}
\bigskip\bigskip
{\bff 1. The question (until 2011)}
\medskip
Let $\,{\cal F}\,$ be the family of all real monotonic functions $\,f\,$
defined on an arbitrary (non\-degenerate) interval $\,I\,$ 
such that $\,f\,$ is {\it singular}, i.e.~$\;f'(x)=0\;$ 
for almost all $\;x\in I\,$.
Let $\,{\cal F}^*\,$ be the family of all functions $\,f\,$ in $\,{\cal F}\,$
such that $\,f'(x)\not=0\,$ for at least one point $\,x\,$
at which $\,f\,$ is differentiable (with a finite derivative).
The classical Cantor function (the {\it devil's staircase})
lies in $\;{\cal F}\setminus {\cal F}^*\,$.
Also the famous    
Riesz-Nagy function (see [1] 18.8)
and Minkowski's {\it Fragefunktion} (see [3] p.~345)
and the interesting function $\,F_{3,2}\,$ investigated 
by Parad¡s e.a.~[3],
which all are strictly increasing and singular, 
have the property that at each point the derivative is 
$\,0\,$ or $\,\infty\,$ or not existent. 
Another classic example of a singular function is given by 
$$F_\varphi(x)\;\,=\;\sum_{\varphi(n)< x}{1\over 2^n}\qquad (x\in{\Bbb R})$$
where $\,\varphi\,$ is any bijection from $\,{\Bbb N}\,$ onto $\,{\Bbb Q}\,$.
(The summation is extended over all $\;n\in{\Bbb N}\;$
with $\;\varphi(n)<x\,$.) The function $\,F_\varphi\,$ 
is the prototype of a strictly increasing function 
which is discontinuous at each rational number
and continuous at each irrational number.
Since obviously $\,F_\varphi\,$ is the limit of
a series of monotonic step functions, we always 
have $\;F_\varphi\in{\cal F}\,$. But again, independently of $\,\varphi\,$, 
we have $\,F_\varphi\not\in{\cal F}^*\,$. 
(See the proof in the appendix.)
\smallskip
Consequently, there arises 
the question arises whether $\;{\cal F}^*=\,\emptyset\,$.
\bigskip\medskip
{\bff 2. The answer (published in 2011)}
\medskip
Let $\,\Phi\,$ be the family of all bijective functions
from $\,{\Bbb N}\,$ onto $\,{\Bbb Q}\,$. 
We modify the definition of $\,F_\varphi\,$ and 
consider functions $\;G_\varphi\,:\;{\Bbb R}\to\,{\Bbb R}\;$ 
for denumerations $\;\varphi\in\Phi\;$
which are defined by
$$G_\varphi(x)\,\;=\;\sum_{\varphi(n)< x}{1\over n^2}\;\,.$$
Of course, just as the functions $\,F_\varphi\,$, 
all functions $\,G_\varphi\,$ 
are strictly increasing and singular (and
continuous precisely at the irrational numbers).
Now, for some $\,\varphi\in\Phi\,$ we actually have 
$\;G_\varphi\in{\cal F}^*\;$ and hence we can be sure that 
$\;{\cal F}^*\not=\,\emptyset\,$. Moreover, we can prove the following
\medskip\smallskip
{\bf Theorem.} $\;$ {\it For every sequence 
of distinct irrational numbers $\;\xi_1,\xi_2,\xi_3,...\,$ 
and every sequence $\;c_1,c_2,c_3,...\;$ of positive real numbers
there is a denumeration $\;\varphi\in\Phi\;$
such that $\,G_\varphi\,$ is differentiable at $\,\xi_k\,$ and
$\;G_\varphi'(\xi_k)=c_k\;$ for every $\;k\in{\Bbb N}\,$.}
\bigskip
Clearly, if $\,A\subset{\Bbb R}\,$ is countable then 
the translate $\;\theta+A\;$ is disjoint from $\,{\Bbb Q}\,$ for some 
$\,\theta\in{\Bbb R}\,$. Therefore, an immediate consequence of the theorem 
is the following 
\medskip
{\bf Corollary.} {\it For every countable set $\,A\subset{\Bbb R}\,$
and every mapping $\,g\,$ from $\,A\,$ into $\;]0,\infty[\;$
there exists a strictly increasing real function $\,f\,$
such that $\;f'(x)=g(x)\;$ for every $\,x\in A\,$ and 
$\;f'(x)=0\;$ for almost all $\,x\in{\Bbb R}\,$.}
\medskip
For example, there exists a strictly increasing real function $\,f\,$
with $\;f'(x)=2^x\;$ for every rational $\,x\,$ and 
$\;f'(x)=0\;$ for almost all irrational $\,x\,$.
(Of course, in this statement the word {\it almost} cannot be omitted.)
\vfill\eject

{\bff 3. The proof}
\medskip
In order to prove the theorem 
define three sequences $\;(a_k),(b_k),(d_k)\;$ of positive 
even numbers
such that with $\;A_k\,=\,\{\,a_k+nd_k\;|\;n\in\Bbb N\,\}\;$ and
$\;B_k\,=\,\{\,b_k+nd_k\;|\;n\in\Bbb N\,\}\;$ all the elements in the family 
$\;\{\,A_k\;|\;k\in\Bbb N\,\}\cup\{\,B_k\;|\;k\in\Bbb N\,\}\;$
are mutually disjoint sets of even numbers. 
(Choose for example 
$\;a_k=2\cdot 3^k\;$ and
$\;b_k=4\cdot 3^k\;$ and
$\;d_k=2\cdot 3^{k+1}\,$.)
For each $\;k\in \Bbb N\;$ define a strictly decreasing sequence
$\,(x_m^{(k)})\,$ which tends to 
0 as $\,m\to\infty\,$ by
$\;x_m^{(k)}\,:=\,(c_kd_k^2m)^{-1}\,.$
Elementary asymptotic analysis yields
$$\sum_{n=m}^\infty {1\over (s+nd_k)^2}\;=\;
{1\over d_k^2m}\,+\,O\Big({1\over m^{2}}\Big)\;=\;
c_k\cdot x_m^{(k)}\,+\,
O\big((x_m^{(k)})^2\big)
\;\;\;(m\to\infty)\leqno(3.1)$$
for each $\;k,s\in\Bbb N\,$.
Now put $\;\delta_1=1\;$ and $\;\delta_k\,:
=\,\min\,\{\,|\xi_i-\xi_k|\;|\;i<k\,\}\;$ for all integers $\,k\geq 2\,$.
Then choose $\,m_k\in\Bbb N\,$ for every $\,k\in\Bbb N\,$ 
such that with $\;y_k\,:=\,x_{m_k}^{(k)}\,=\,(c_kd_k^2m_k)^{-1}\;$
\medskip
(3.2)\qquad\qquad\qquad\qquad\qquad\qquad\qquad
$y_k+\sqrt{c_ky_k} \;<\;\delta_k/2$
\medskip
(3.3)\qquad\qquad\qquad\qquad\qquad\qquad
\qquad$2\,c_{k+1}\,y_{k+1}\;\leq\;c_k\,y_k$
$$\max\Big\{\,\sum_{n=m_k}^\infty {1\over (a_k+nd_k)^2}\,,\;
\sum_{n=m_k}^\infty {1\over (b_k+nd_k)^2}\,\Big\}\,\;\leq\,\;
2\,c_k\,y_k\leqno(3.4)$$
for each $\,k\in\Bbb N\,$.
Define sets
$\,X_k\subset A_k\,$ and $\,Y_k\subset B_k\,$
by $\;X_k\,=\,\{\,a_k+nd_k\;|\;m_k\leq n\in\Bbb N\,\}\;$ and
$\;Y_k\,=\,\{\,b_k+nd_k\;|\;m_k\leq n\in\Bbb N\,\}\,$.
Now we define our desired denumeration $\,\varphi\in\Phi\,$ 
firstly on the domain 
$\;D\,:=\,\bigcup\{\,X_k\cup Y_k\;|\;k\in\Bbb N\,\}\;$
by choosing for each $\;k\in\Bbb N\;$ and every integer $\;n\geq m_k\;$
\smallskip 
\centerline{$\varphi(a_k+nd_k)\;\in\;\,
]\xi_k+x_{n+1}^{(k)},\xi_k+x_{n}^{(k)}[\;\cap\,\Bbb Q\setminus \Bbb Z$}

and

\centerline{$\varphi(b_k+nd_k)\;\in\;\,
]\xi_k-x_{n}^{(k)},\xi_k-x_{n+1}^{(k)}[\;\cap\,\Bbb Q\setminus \Bbb Z\,$.}
\medskip 
Of course,
these choices can be made so that $\,\varphi\,$ is injective on 
the domain $\,D\,$.
(Choose for example for each $\,k\,$ only rationals of the form 
$\,r\cdot p_k^s\,$ with $\,r,s\in{\Bbb Z}\,$ where $\,p_k\,$
is the $k$th prime number.) 
Then we extend $\,\varphi\,$ by defining $\,\varphi^{-1}\,$
on $\;\Bbb Q\setminus(\Bbb Z\cup\varphi(D))\;$ 
via 
\smallskip
\centerline{$\varphi^{-1}\big({p\over q}\big)\;:=\;
 {\sqrt{13}}^{1+|p|/p}\cdot 3^{|p|}\cdot 5^q\cdot 7^{\delta(p/q)}$}
\smallskip
where $\,p,q\,$ are coprime integers and $\,q\geq 2\,$ and where
$\,\delta (p/q)\,$ is the least positive integer not smaller than 
$\;\max\,\{\,|{p\over q}-\xi_i|^{-1}\;\,|\;\,i=1,2,...,q\,\}\,$.
This extension is clearly possible because
$\,\varphi^{-1}\big({p\over q}\big)\,$ is always odd 
by definition 
and $\,D\,$ contains only even numbers.
(Using the primes 3, 5, 7, 13 we can be sure that $\,\varphi\,$
is well-defined and injective.)
Finally we extend $\,\varphi\,$ in any way to a bijection
from $\,\Bbb N\,$ onto $\,\Bbb Q\,$.
\medskip
Now fix $\;\kappa\in\Bbb N\;$ and for abbreviation put
$\;\xi:=\xi_\kappa\;$ and
$\;x_m:=x_m^{(\kappa)}\;$ for every $\,m\in\Bbb N\,$.
For $\;k,m\in\Bbb N\;$ let 
$\;\,{\cal I}_{m,k}\,:=\,[\xi_k-y_k,\xi_k+y_k]\,\cap\,[\xi-x_m,\xi+x_m]\;$.
We claim that
\smallskip
(3.5)\qquad $\forall\,m,k\in{\Bbb N}\,:\;
k>\kappa\; \land \;{\cal I}_{m,k}\not=\emptyset\;\Longrightarrow\;
c_ky_k\leq x_m^2\;$.
\smallskip
Indeed, if $\;k>\kappa\;$ then $\;\xi\not\in [\xi_k-y_k,\xi_k+y_k]\;$
since by (3.2) $\;y_k<\delta_k\leq|\xi_k-\xi|\,$.
Therefore, if additionally $\;{\cal I}_{m,k}\not=\emptyset\;$ for any $\,m\,$
then we clearly must have $\;x_m+y_k\,\geq\,|\xi_k-\xi|\,\geq\,\delta_k\;$
and thus $\;c_ky_k>x_m^2\;$ would imply
$\;\sqrt{c_ky_k}+y_k>x_m+y_k\geq\delta_k\;$ contrarily to (3.2).
\medskip
In order to conclude the proof by verifying 
$\;G_\varphi'(\xi)=c_\kappa\;$ we 
take into account the following three considerations. 
\vfill\eject
\medskip
Firstly, there clearly exists a bound
$\,\delta>0\,$ such that 
$\;\Bbb Z\,\cup\,\varphi\big(\bigcup_{k=1}^{\kappa-1}(X_k\cup Y_k)\big)\;$
is disjoint from $\;[\xi-\delta,\xi+\delta]\,$. 
We claim that the set 
\smallskip
\centerline{$\;K_m\,:=\,\{\,k\in\Bbb N\;\,|\,\;k>\kappa\;\,\land\,\;
{\cal I}_{m,k}\not=\emptyset\,\}\;$} 
\smallskip
is empty for some $\;m\in\Bbb N\;$
if and only if $\,\xi\,$ is not a limit point of the set 
$\;\{\xi_1,\xi_2,\xi_3,...\}\,$.                         
Indeed, if $\,K_m\,$ is empty for some $\,m\,$ then for every $\,k>\kappa\,$
we have $\,{\cal I}_{m,k}=\emptyset\,$ and hence 
$\;\xi_k\not\in[\xi-x_m,\xi+x_m]\,$, whence
$\,\xi\,$ cannot be a limit point of $\;\{\xi_1,\xi_2,\xi_3,...\}\,$. 
Conversely, if $\,\xi\,$ is not a limit point then 
we may choose $\,h>0\,$ so that $\;\xi_k\not\in[\xi-h,\xi+h]\;$
for every $\,k>\kappa\,$. Since by (3.2) we have 
$\;y_k<{1\over 2}|\xi_k-\xi|\;$ for every 
$\;k>\kappa\,$, we must have $\,{\cal I}_{m,k}=\emptyset\,$ 
for every $\,k>\kappa\,$ or, equivalently, $\,K_m=\emptyset\,$
if $\,m\,$ is chosen so that $\,x_m<{h\over 2}\,$.
\smallskip
So if $\,\xi\,$ is not a limit point of the set 
$\;\{\xi_1,\xi_2,\xi_3,...\}\;$ then there exists a number $\,\tilde m\,$
such that $\;{\cal I}_{\tilde m,k}=\emptyset\;$ for every $\,k>\kappa\,$ 
and hence $\;\tilde \delta\,=\,\min\{\delta, x_{\tilde  m}\}\;$ 
is a bound such that even 
$\;\Bbb Z\,\cup\,\varphi\big(\bigcup_{k\not=\kappa}(X_k\cup Y_k)\big)\;$
is disjoint from $\;[\xi-\tilde\delta,\xi+\tilde\delta]\,$.
\medskip
Secondly, assume that $\,\xi\,$ is a limit point of 
$\;\{\xi_1,\xi_2,\xi_3,...\}\,$. Then 
$\,K_m\,$ is never empty and we may define 
$\;\mu(m)\,:=\,\min\,K_m\,$. 
Put 
\smallskip
\centerline{$\;L_{m,k}\,:=\,
\{\,n\in (X_k\cup Y_k)\;\,|\,\;\xi-x_m\leq \varphi(n) < \xi+x_m\,\}\,$.}
\smallskip
Then in view of the definition of $\,\varphi\,$, 
$$\sum\limits_{k>\kappa}\sum_{n\in L_{m,k}}{1\over n^2}
\;\leq\;\sum\limits_{k\in K_m}\Big(\sum_{n=m_k}^\infty{1\over (a_k+nd_k)^2}
\,+\,\sum_{n=m_k}^\infty{1\over (b_k+nd_k)^2}\Big)$$
for all $\,m\in{\Bbb N}\,$. Thus by applying (3.4),
$$\sum\limits_{k>\kappa}\sum_{n\in L_{m,k}}{1\over n^2}
\;\leq\;4\sum\limits_{k\in K_m}c_ky_k\;.$$
Furthermore, since 
$\;c_{\mu(m)+n}y_{\mu(m)+n}\leq 2^{-n} c_{\mu(m)}y_{\mu(m)}\;$
for $\;n\,=\,0,1,2,3,...\;$ due to (3.3), 
$$\sum\limits_{k\in K_m}c_ky_k\;\leq\;
\sum\limits_{k=\mu(m)}^\infty c_ky_k\;\leq\;
\sum\limits_{n=0}^\infty 2^{-n}c_{\mu(m)}y_{\mu(m)}
\;=\;2c_{\mu(m)}y_{\mu(m)}\;.$$
By (3.5) we have $\;c_{\mu(m)}y_{\mu(m)}\;\leq\;x_m^2\;$
and so putting all together we arrive at  
$${1\over x_m}\cdot
\sum\limits_{k>\kappa}\sum_{n\in L_{m,k}}{1\over n^2}
\;\leq\;8x_m\;\to\;0\quad(m\to\infty)\leqno(3.6)$$
provided that $\,\xi\,$ is a limit point of 
$\;\{\xi_1,\xi_2,\xi_3,...\}\,$.
\medskip
Thirdly, define
$\;N_m\,:=\,\{\,n\in\Bbb N\setminus D\;\,|\;\,
\xi-x_m\leq \varphi(n) < \xi+x_m\,\}\;$
and let $\,M\,$ be the smallest positive integer such that
the interval $\,[\xi-x_M,\xi+x_M]\,$ 
does not contain
integers or reduced fractions $\,{p\over q}\,$ with $\,|q|<\kappa\,$. 
Then for each $\,m\geq M\,$ we have 
$\;N_m\subset [5\cdot 7^{1/x_m},\infty[\;$
because if $\,n\in N_m\,$ and $\,\varphi(n)={p\over q}\,$
(where $\,q>0\,$ and the fraction $\,{p\over q}\,$ is reduced)
then $\;{p\over q}\in[\xi-x_m,\xi+x_m]\subset [\xi-x_M,\xi+x_M]\;$
and hence (by the definition of $\,M\,$) $\,\kappa\in\{1,2,...,q\}\,$
and hence $\;{1\over x_m}\leq |{p\over q}-\xi|^{-1}\leq \delta(p/q)\;$
and hence $\;n\,\geq\,3^{|p|}\cdot 5^q\cdot 7^{\delta(p/q)}
\,\geq\,5\cdot 7^{1/x_m}\,$.
Consequently, for $\;m\geq M\;$
$${1\over x_m}\cdot
\sum\limits_{n\in N_m} {1\over n^2}\;\leq\;
{1\over x_m}\cdot\!\!\!
\sum\limits_{n\geq 5\cdot 7^{1/x_m}}\!\!{1\over n^2}\;\leq\;
{1\over x_m}\cdot\!\!\!\int\limits_{4\cdot 7^{1/x_m}}^\infty
\!\!{{\rm d}x\over x^2}
\;\,\to\,\;0\;\;\;(m\to\infty)\;.\leqno(3.7)$$
\medskip\smallskip
In order to conclude the proof by verifying 
$\;G_\varphi'(\xi)=c_\kappa\;$ 
it is enough to verify 
$$\lim\limits_{m\to\infty}
{G_\varphi(\xi+x_m)-G_\varphi(\xi)\over x_m}
\;\,=\,\;
\lim\limits_{m\to\infty}
{G_\varphi(\xi)-G_\varphi(\xi-x_m)\over x_m}\;\,=\,\;c_\kappa$$
because $\;\lim\limits_{m\to\infty}{x_{m+1}\over x_m}=1\;$
and $\,G_\varphi\,$ is increasing.
\medskip
Now, for every $\,m\in{\Bbb N}\,$ we can write 
$${G_\varphi(\xi+x_m)-G_\varphi(\xi)\over x_m}\;=\;
{1\over x_m}\cdot\Big(\sum\limits_{n\in S_m\cap X_\kappa}{1\over n^2}\,+\,
\sum\limits_{n\in S_m\cap D\setminus X_\kappa}{1\over n^2}\,+\,
\sum\limits_{n\in S_m\setminus D}{1\over n^2}\;\Big)$$
where $\;S_m\,:=\,\{\,n\in{\Bbb N}\;|\;\xi\leq \varphi(n)<\xi+x_m\,\}\,.$
(Recall that $\,X_\kappa\subset D\,$.)
In view of $\;S_m\setminus D\,\subset\,N_m\;$ and (3.7) we have
$$\lim\limits_{m\to\infty}\;                                    
{1\over x_m}\;\cdot\sum\limits_{n\in S_m\setminus D}{1\over n^2}
\,\;=\,\;0\,\,.$$
In view of (3.1) and the definition of $\,\varphi\,$ 
we have 
$$\lim\limits_{m\to\infty}\;
{1\over x_m}\;\cdot\sum\limits_{n\in S_m\cap X_\kappa}{1\over n^2}
\,\;=\;\,c_\kappa\;.$$
Since $\;S_m\cap D\setminus X_\kappa\;$ is a subset of 
$\,\bigcup_{k\not=\kappa} (X_k\cup Y_k)\,$, we have
$$\lim\limits_{m\to\infty}\;                                    
{1\over x_m}\;\cdot
\sum\limits_{n\in S_m\cap D\setminus X_\kappa}{1\over n^2}\,\;=\,\;0$$
in view of (3.6) and the 
consideration involving the bound $\,\delta\,$ 
and the potential bound $\,\tilde\delta\,$.
(Clearly, if
$\,\xi\,$ is not a limit point of $\;\{\xi_1,\xi_2,\xi_3,...\}\;$ 
then $\;S_m\cap D\setminus X_\kappa\,=\,\emptyset\;$ 
for sufficiently large $\,m\,$.)
Summing up,
$$\lim\limits_{m\to\infty}
{G_\varphi(\xi+x_m)-G_\varphi(\xi)\over x_m}\;\,=\,\;c_\kappa\;.$$
Analogously,
$$\lim\limits_{m\to\infty}
{G_\varphi(\xi)-G_\varphi(\xi-x_m)\over x_m}\;\,=\,\;c_\kappa$$
and this finishes the proof.
\bigskip\bigskip\medskip
{\it Remark.} It is not true that for
every $\,\varphi\in\Phi\,$ there are points $\,\xi\,$ 
such that $\,0<G_\varphi'(\xi)<\infty\,$. 
(Choose any $\,\varphi\in\Phi\,$ where
$\,\varphi\,$ maps $\;\{\,2^n\;|\;n\in\Bbb N\,\}\;$
onto $\;\Bbb Q\setminus \Bbb Z\;$ and define $\,\psi\in\Phi\,$ anyhow
so that $\,\psi(m)=\varphi(2^{m/2})\,$ for every even $\,m\in{\Bbb N}\,$.
Then for every $\,k\in{\Bbb Z}\,$ there is a constant $\,\tau_k\,$
such that $\;G_\varphi(x)=F_\psi(x)+\tau_k\;$ whenever $\,k<x\leq k+1\,$.
Consequently, by the proposition below, $\,G_\varphi\not\in {\cal F}^*\,$.)
\vfill\eject
{\bff 4. Appendix}
\medskip
In the following we prove the statement from the first chapter. 
\medskip
{\bf Proposition.} {\it Independently 
of the denumeration $\;\varphi\in\Phi\,$,
there never exists 
a real $\,\xi\,$ such that $\,F_\varphi\,$ is differentiable at $\,\xi\,$
and $\;F_\varphi'(\xi)\not=0\,$.}
\medskip
{\it Proof.} Since $\,F_\varphi\,$ is increasing,
$\;F_\varphi'(\xi)\not=0\;$ means 
$\;F_\varphi'(\xi)>0\,$.  
Suppose that $\;F_\varphi'(\xi)=2^x\;$ for an irrational 
$\,\xi\,$ and a real $\,x\,$.
Fix $\;0<\varepsilon<{1\over 10}\;$ and 
$\;N\in{\Bbb N}\;$ such that 
$$2^{x-\varepsilon}\;<\;
{F_\varphi(\xi\pm 2^{-m})-F_\varphi(\xi)\over \pm 2^{-m}}\;<\; 
2^{x+\varepsilon}$$
for every $\;m\geq N\,$.
For each $\;m\geq N\;$ define
\smallskip
\centerline{${\cal N}_+(m)\;:=\;\{\,n\in{\Bbb N}\;\,|\,\;\xi\leq \varphi(n)<\xi+2^{-m}\,\}\;,$}
\smallskip
\centerline{${\cal N}_-(m)\;:=\;\{\,n\in{\Bbb N}\;\,|\,\;\xi-2^{-m}\leq \varphi(n)<\xi\,\}\;,$}
\smallskip
whence
$$2^{x-\varepsilon}\;<\;
2^m\sum_{n\in{\cal N}_\pm(m)}2^{-n}\;<\; 2^{x+\varepsilon}$$
for all $\;m\geq N\,$.
If $\,{\cal N}\,$ is a nonempty subset of $\,{\Bbb N}\,$ with 
minimum $\,\mu\,$ then of course
\smallskip
\centerline{$2^{-\mu}\;\leq\;\sum\limits_{n\in{\cal N}}2^{-n}\;\leq\;2^{1-\mu}\;$.}
\smallskip
Thus with $\;\mu_{\pm}(m)\,:=\,\min{\cal N}_\pm(m)\;$ we must have 
$\;m-x-\varepsilon\,<\,\mu_{\pm}(m)\,<\,1+m-x+\varepsilon\;$ 
for all $\;m\geq N\,$.
Therefore we must have
\smallskip
\centerline{$[m-x-\varepsilon\,,\,1+m-x+\varepsilon]\,\cap\,
({\cal N}_+(m)\cup{\cal N}_-(m))\;=\;
[m-x-\varepsilon\,,\,1+m-x+\varepsilon]\,\cap\,{\Bbb Z}\;$}
\smallskip
for all $\;m\geq N\;$
because the interval $\;[m-x-\varepsilon\,,\,1+m-x+\varepsilon]\;$
contains at most two integers and   $\;\mu_+(m)\not=\mu_-(m)\;$
for all $\;m\geq N\,$. (Naturally,
the two sets  $\,{\cal N}_+(m)\,$ and
$\,{\cal N}_-(m)\,$ are always disjoint.)
Now obviously
\smallskip
\centerline{$\;\bigcup\limits_{m\geq N}[m-x-\varepsilon\,,\,1+m-x+\varepsilon]\cap{\Bbb Z}
\;\,=\,\;[N-x-\varepsilon\,,\,\infty[\;\cap\,{\Bbb Z}\,.$}
\smallskip
Consequently, 
$\;{\cal N}_+(N)\cup{\cal N}_-(N)\;\supset
\;[N-x-\varepsilon\,,\,\infty[\;\cap\,{\Bbb Z}\;$
since $\;{\cal N}_\pm(m)\subset {\cal N}_\pm(N)\;$ 
for every $\;m\geq N\,$. 
Therefore the set 
$\;{\Bbb N}\,\setminus\,({\cal N}_+(N)\cup{\cal N}_-(N))\;$ must be finite,
but this is impossible since there
are infinitely many rationals outside the interval 
$\;[\xi-2^{-N},\xi+2^{-N}]\;$ which have to be numbered
by $\;\varphi\,$, {\it q.e.d.}
\bigskip\bigskip
{\it Remark.} In view of the preceding proof it is plain
that there is at most one point $\,\xi\,$ at which
the right or left derivative of $\,F_\varphi\,$
exists and is not equal to $\,0\,$ or $\,\infty\,$.
(Note that the left derivative may exist at $\,\xi\,$
even when $\,\xi\in{\Bbb Q}\,$.)
Moreover, by choosing $\,\varepsilon\,$ so that 
$\;[x-\varepsilon,x+\varepsilon]\cap{{\Bbb Z}}\,=\,\emptyset\,$,
it is easy to verify that the 
right or left derivative of $\,F_\varphi\,$
can never assume a value 
$\,2^x\,$ with $\;x\not\in{{\Bbb Z}}\,$.
But a more detailed investigation (see [2]) shows that 
values $\,2^x\,$ with $\;x\in{{\Bbb Z}}\;$
cannot be assumed either. 
Thus the proposition can be generalized so that {\it neither the 
right nor the left derivative of $\,F_\varphi\,$ is finite and positive
anywhere.}
\vfill\eject
{\bff References}
\bigskip
[1] E.H.~{\smc Hewitt} and K.~{\smc Stromberg}:
Real and Abstract Analysis. Springer, 1965.
\medskip
[2] G.~{\smc Kuba}: {\it On the differentiability of certain saltus 
functions}. 

\rightline{Colloq.~Math.~{\bf 125}, 15-30 (2011).}  
\smallskip
[3] J.~{\smc Parad¡s}, P.~{\smc Viader}, and L.~{\smc Bibiloni}:
{\it A New Singular Function.}

\rightline{Am.~Math.~Monthly {\bf 118}, 344-354 (2011).}

\bigskip\bigskip

{\bf Author's address:} Gerald Kuba, Institute of Mathematics,  

University of Natural Resources and Life Sciences, Vienna, Austria
\smallskip
E-mail: {\tt gerald.kuba@boku.ac.at}

\bigskip\medskip\bigskip\medskip
\hrule
\bigskip\medskip\bigskip\medskip
{\bf The theorem and its proof is contained 
in the author's paper [2].}

\end

{\bff References}